\documentclass[11pt,letterpaper]{amsart}

\usepackage{tikz}
\tikzstyle{every node}=[circle, draw, fill=black!50,
                        inner sep=0pt, minimum width=3pt]

\usepackage{amsfonts,amsmath,latexsym,color,epsfig,hyperref,enumitem, amssymb,bbm}

\newtheorem{theorem}{Theorem}[section]
\newtheorem{lemma}{Lemma}[section]

\newtheorem{prop}[lemma]{Proposition}

\renewcommand{\le}{\leqslant}
\renewcommand{\ge}{\geqslant}

\def\qed{\ifvmode\mbox{ }\else\unskip\fi\hskip 1em plus 10fill$\Box$}

\input{epsf}
\makeatletter
\def\Ddots{\mathinner{\mkern1mu\raise\p@
\vbox{\kern7\p@\hbox{.}}\mkern2mu
\raise4\p@\hbox{.}\mkern2mu\raise7\p@\hbox{.}\mkern1mu}}
\makeatother

\def\R{\mathbb R}

\def\E{\mathbb E}

\def\IS{\int_{S^{n-1}} }

\title{ Spherical sets avoiding orthonormal bases}
\author{Dmitrii Zakharov}
\thanks{Zakharov's research was supported by the Jane Street Graduate Fellowship.}
\address{Department of Mathematics, Massachusetts Institute of Technology, Cambridge, MA 02139, USA}
\email{zakhdm@mit.edu}

\date{}

\begin{document}
\maketitle

\begin{abstract}
    We show that there exists an absolute constant $c_0<1$ such that for all $n \ge 2$, any measurable set $A \subset S^{n-1}$ of density at least $c_0$ contains $n$ pairwise orthogonal vectors. The result is sharp up to the value of the constant $c_0$. 

    Moreover, we show that for all $2\le k \le n$ a set $A$ avoiding $k$ pairwise orthogonal vectors has measure at most $\exp(-c_1 \min\{\sqrt{n}, n/k\})$ for some $c_1>0$.
    Proofs rely on the harmonic analysis on the sphere and the hypercontractive inequality.
\end{abstract}

\section{Introduction}

We are interested in the following general question: given $n\ge 2$ and a finite set of points $P$ in $S^{n-1}$ or $\R^n$, what is the largest density of a subset $A$ in $S^{n-1}$, resp. in $\R^n$, which does not contain a congruent copy of the set $P$?\footnote{Here and throughout all sets are assumed to be measurable.} The most classical and well-studied case of this question is when $P = \{x, y\}$ is a 2-element set of points. Then one wants to find the largest possible density $m_{S^{n-1}}(t)$, resp. $m_{\R^n}(t)$, of a set $A$ with a forbidden distance $t = \operatorname{dist}(x,y)$. In the Euclidean space $\R^n$, all distances are equivalent, so it is enough to consider $t=1$. The celebrated result of Frankl--Wilson \cite{Frankl}, later improved upon in \cite{Rai2}, implies that $m_{\R^n}(1)$ decreases exponentially in $n$. In the plane $\R^2$, 
it was recently shown \cite{Ambrus} that any set avoiding the unit distance has density at most $0.25-\varepsilon$ where $\varepsilon =0.003 > 0$. This answered a question of Erd{\H o}s and gave a quantitative improvement of the fact that the measurable chromatic number of $\R^2$ is at least 5, i.e. that one cannot partition $\R^2$ into four measurable sets avoiding unit distances (earlier, De Grey \cite{DeGrey} showed that this holds even without the measurability assumption).

More generally, there has been a lot of effort \cite{Bachoc, Castro2, DeCorte, Oliv} to obtain better bounds on these densities in small dimensions. The techniques developed in this line of work rely heavily on harmonic and Fourier analysis and linear and semidefinite programming.

On the sphere $S^{n-1}$, the choice of the distance $t$ is important. For the perhaps most natural choice $t = \sqrt{2}$, i.e. when we forbid our set to contain pairs of orthogonal vectors, Kalai's double cap conjecture \cite{kalai} predicts that the largest set with this property is
$$
A_{0} = \left\{x \in S^{n-1}: |x_1| > \frac{1}{\sqrt{2}}\right\}.
$$
For $n=3$, one has $|A_0| \approx 0.292$ and the best current upper bound \cite{Bekker} is roughly $0.297$. For large $n$, the measure of $A_0$ is asymptotically $(2+o(1))^{-n/2}$ and the Frankl--Wilson's method implies an exponential upper bound on $m_{S^{n-1}}(\sqrt{2})$. In fact, Raigorodskii \cite{Rai} showed that for any fixed $t \in (0,2)$, the function $m_{S^{n-1}}(t)$ decays exponentially in $n$.

Less is known when the forbidden set $P$ has size greater than $2$. Let $m_{S^{n-1}}(P)$ and $m_{\R^n}(P)$, denote the maximal density of a set $A$ in $S^{n-1}$, resp. $\R^n$, avoiding a congruent copy of $P$. In \cite{Castro2}, Castro-Silva, de Oliveira, Slot and Vallentin introduce a semidefinite programming approach to this `higher uniformity' problem. In large dimensions, their approach implies the following. For $k\ge 2$ and $t \in (-1,1)$, let $\Delta_{k,t}$ denote the set of $k$ unit vectors with pairwise scalar product $t$, note the switch from distance to scalar product which is a more natural quantity on the sphere. Then for any fixed $k \ge 2$ and $t \in (0,1)$, $m_{S^{n-1}}(\Delta_{k,t})$ decays exponentially in $n$. They also showed that $m_{\R^n}(\Delta_{k,0}) \le (1-\frac{1}{9k^2}+o(1))^n$ for any fixed $k$; this bound was later improved by combinatorial methods in \cite{KSZ}.

Neither of the approaches in \cite{Castro2} nor \cite{KSZ} lead to particularly strong bounds on $m_{S^{n-1}}(\Delta_{k,t})$ in the case when the forbidden scalar product $t$ is non-positive or the number of forbidden points $k$ is close to $n$. For instance, the recursive semidefinite bound from \cite{Castro2} only leads to the trivial upper bound $m_{S^{n-1}}(\Delta_{k,0}) \le \frac{k-1}{n}$. To see this, let $X = \{x_1, \ldots, x_n\} \subset S^{n-1}$ be a uniformly random collection of $n$ pairwise orthogonal vectors on $S^{n-1}$ (which can be obtained e.g. by applying a uniformly random rigid motion to the standard basis). Then by linearity of expectation, we have $\E |A \cap X| = n \mu(A)$. On the other hand, if $A$ has no $k$ pairwise orthogonal vectors then $|A\cap X| \le k-1$ for any $X$, and so $\E|A\cap X| \le k-1$ holds, giving the desired bound.
In particular, when $k = n$ this gives an upper bound of $|A| \le 1- \frac{1}{n}$ for a set $A\subset S^{n-1}$ which does not contain $n$ pairwise orthogonal vectors. On the other hand, we have the following construction:
\begin{equation}\label{A1}
A_1 = \left\{x \in S^{n-1}:~ |x_1| < \frac{1}{\sqrt{n}}\right\}.    
\end{equation}
We claim that $A_1$ does not contain $n$ pairwise orthogonal vectors.
Indeed, suppose that $v_1, \ldots, v_n \in A_1$ are pairwise orthogonal. But if we let $e_1 = (1, 0, \ldots, 0)$ then:
$$
1 = \|e_1\|^2 = \langle e_1, v_1\rangle^2 + \ldots + \langle e_1, v_n\rangle^2 < n \left(\frac{1}{\sqrt{n}}\right)^2 = 1,
$$
a contradiction. For large $n$, the measure of the set $A_1$ approaches $\frac{1}{\sqrt{2\pi}}\int_{-1}^1 e^{-t^2/2}dt \approx 0.68$, which is the probability that a Gaussian random variable does not exceed its variance. 
In this note, we address the limitations of the previous approaches and show that, up to a constant, the example (\ref{A1}) is best possible and that the simple upper bound $1-\frac{1}{n}$ is very far from the truth.

\begin{theorem}\label{main}
There exists an absolute constant $c_0 < 1$ such that for any $n \ge 2$, any subset $A \subset S^{n-1}$ of measure at least $c_0$ contains $n$ pairwise orthogonal vectors.  
\end{theorem}

The set $A_1$ above demonstrates that one cannot take $c_0 < 0.68$ in Theorem \ref{main}. Our proof goes through with something like $c_0 = 1-10^{-6}$ but we did not attempt to optimize this value.

Using similar approach we can show upper bounds on $m_{S^{n-1}}(P)$ for various patterns $P$ on the sphere whose vectors are linearly independent and the pairwise scalar products are close to zero. For simplicity we only give the bound in the case of pairwise orthogonal vectors:

\begin{theorem}\label{prop:other_simplices}
    There is a constant $c_1 >0$ such that for any $2\le k\le n$ we have 
    \[
    m_{S^{n-1}}(\Delta_{k,0}) \le \exp\left(-c_1\min\left\{\sqrt{n}, \frac{n}{k}\right\}\right).
    \]
\end{theorem}

Note that taking $k=n$ recovers Theorem \ref{main}. 
On the other hand, the natural `double cap' construction gives a lower bound 
\[
m_{S^{n-1}}(\Delta_{k,0}) \ge \left(1 - \frac{1}{k}+o(1)\right)^{n/2},
\]
where $o(1)$ tends to 0 with $k/n \rightarrow 0$. So for $k \gtrsim \sqrt{n}$ we get matching behavior but for $k \lesssim \sqrt{n}$ there is a gap. 

\subsection{Hyperplane slices.} Theorem \ref{prop:other_simplices} states that we can find a large collection of pairwise orthogonal vectors in a set $A \subset S^{n-1}$ of an appropriate density. We construct these vectors inductively: we start by picking a vector $x_1 \in A$, then we pick $x_2 \in A \cap x_1^\perp$, and then we pick $x_3 \in A \cap x_1^\perp \cap x_2^\perp$ and so on (here and throughout $x^\perp$ stands for the hyperplane orthogonal to a vector $x$). If we can ensure that at each step the intersection $A \cap x_1^\perp \cap \ldots \cap x_j^\perp$ is non-empty, then after $k-1$ steps we will produce $k$ pairwise orthogonal vectors in $A$.

In order to execute this strategy, we show that if $A$ is sufficiently dense then we can find $x_1 \in A$ so that the density the set $A \cap x_1^{\perp}$ (as a subset of an $(n-2)$-dimensional sphere) has a good lower bound. In what follows, we will often use expressions of the form $\mu_{S^{n-2}}(A \cap x^\perp)$ when referring to the density of $A\cap x^\perp$ with respect to the sphere $S^{n-1} \cap x^\perp$.

The following two lemmas are designed to accomplish this in two different ranges of densities. Before stating the lemmas, let first us give a more elementary bound which will not be sufficient for us:

\begin{prop}\label{prop:easy}
    Let $A \subset S^{n-1}$ be a set of measure $\alpha \in (0,1)$, then we have 
    \[
    \E_{x \in A} \mu_{S^{n-2}}(A \cap x^{\perp}) \ge \alpha - \frac{1-\alpha}{n-1}.
    \]
\end{prop}

This Proposition is not new and essentially appears in many linear programming approaches to forbidden configuration problems on $S^{n-1}$, see e.g. this is a simple corollary of the theta-function method in \cite{Bekker}, \cite{Castro2}. We give two proofs of this result: one using the elementary probabilistic approach we used to show that $m_{S^{n-1}}(\Delta_{k,0}) \le \frac{k-1}{n}$ and another one using harmonic analysis on the sphere, which is the main tool of our work and which will be developed in Section \ref{section:prelim}.

\begin{proof}[First proof of Proposition \ref{prop:easy}]
    As before let $X = \{x_1, \ldots, x_n\} \subset S^{n-1}$ be a uniformly random collection of $n$ pairwise orthogonal vectors. Then by linearity of expectation we have $\E |X\cap A| = n \alpha$. On the other hand, we can compute
    \[
    \E {|X\cap A|\choose 2} = \sum_{1\le i < j \le n} \E_X 1_{x_i \in A} 1_{x_j \in A} = {n \choose 2} \alpha \E_{x \in A} \mu_{S^{n-2}}(A \cap x^{\perp}).
    \]
    Here we used the fact that if we condition on $x_i$ then $x_j$ is uniformly distributed on $S^{n-1} \cap x_i^\perp$.
    
    So using convexity of $t \mapsto {t \choose 2}$ we get
    \[
    {n \choose 2} \alpha \E_{x \in A} \mu_{S^{n-2}}(A \cap x^{\perp}) \le {n \alpha \choose 2}
    \]
    which after rearranging gives the desired bound. 
\end{proof}

It turns out that Proposition \ref{prop:easy} is far from being sharp if density of $A$ is close to 1 or 0.

\begin{lemma}\label{lemperp}
    There exist absolute constants $\varepsilon_{\ref{lemperp}} \in (0,1)$ and $C_{\ref{lemperp}} \ge 1$ such that the following holds for all $\varepsilon \le \varepsilon_{\ref{lemperp}}$.
    Let $n \ge 2$ and $A \subset S^{n-1}$ be a centrally symmetric set of density $\alpha = 1-\varepsilon$. Then there exists a point $x \in A$ such that 
    \begin{equation}\label{perp}
        \mu_{S^{n-2}}(A \cap x^\perp) \ge \alpha - \frac{C_{\ref{lemperp}} \varepsilon^{1/2}}{n^2}.
    \end{equation}
\end{lemma}

\begin{lemma}\label{lem:density_lem}
Let $\alpha \in (0, e^{-2})$ and let $A \subset S^{n-1}$ be a set of measure $\alpha$. Then there exists a point $x \in A$ such that for some constant $C_{\ref{lem:density_lem}}>0$:
\[
\mu_{S^{n-2}}(A \cap x^\perp) \ge \alpha \left(1- \frac{C_{\ref{lem:density_lem}} \log(1/\alpha)^2}{n} \right).
\]
\end{lemma}

Let us point out that Lemma \ref{lem:density_lem} is still true if we take the average over all $x \in A$ instead of picking one. On the other hand, this is not the case for Lemma \ref{lemperp}: if we take $A$ to be a spherical band of density $1-\varepsilon$ then the average density of the intersection $A\cap x^\perp$ is of order $\alpha - \Theta(\frac{\varepsilon^2 \log^2 (1/\varepsilon)}{n})$. So a crucial step in the proof of Lemma \ref{lemperp} is to bias the uniform distribution on $A$ in order to increase the average density.

Theorems \ref{main} and \ref{prop:other_simplices} follow from these lemmas by an inductive argument. In case of Theorem \ref{main}, we apply Lemma \ref{lemperp} repeatedly to construct a sequence of $n-n_0$ pairwise orthogonal points in $A$ and then use the trivial bound $\mu_{S^{n_0-1}}(\Delta_{n_0,0}) \le 1-1/n_0$ to construct the remaining $n_0$ points. In case of Theorem \ref{prop:other_simplices}, we use Theorem \ref{main} to reduce to the case $k \le cn$ for a small constant $c>0$ and then iterate Lemma \ref{lem:density_lem}.

To prove the Lemmas, we use harmonic analysis on the sphere and crucially rely on the hypercontractive inequality for on the sphere. Hypercontractivite inequalities are a ubiquitous tool in boolean analysis, in particular, they has been effectively applied to study forbidden intersection problems in discrete product-like spaces, see e.g. \cite{Keevash, Keller, Keller2}.
We hope that there will be more applications of these ideas to continuous forbidden subconfiguration problems as well.

We should point out that similar computations with measures of hyperplane slices on the sphere were also used in other contexts, see for example \cite{klartag2010quantum}.

{\em Acknowledgments.} I thank Mehtaab Sawhney for the help with probabilistic estimates. I thank Danila Cherkashin and anonymous referees for helpful comments.

\section{Preliminaries}\label{section:prelim}

\subsection{Harmonic analysis on the sphere} We recall some basic facts about the space of functions on the sphere, see e.g. \cite{Abram, Castro1, Dai} for a more comprehensive account of the theory.  
Let $L^2(S^{n-1})$ denote the Hilbert space of square-integrable real-valued functions on the $(n-1)$-dimensional unit sphere $S^{n-1} \subset \R^n$. Let $\mu$ be the probability measure on $S^{n-1}$ and define the scalar product of functions $f,g \in L^2(S^{n-1})$ to be
$$
\langle f, g \rangle = \IS f(x) g(x) d\mu(x).
$$
For $d \ge 0$ let the $\mathcal H_{n,d}$ be the space of homogeneous harmonic polynomials of degree $d$ in $n$ variables. That is, a degree $d$ homogeneous polynomial $f \in \R[x_1, \ldots, x_n]$ belongs to $\mathcal H_{n, d}$ if 
$$
\Delta f = \sum_{i=1}^n \frac{\partial^2}{\partial x_i^2} f = 0.
$$
The dimension of this space is given by
$$
\dim \mathcal H_{n, d} = {n+d-1 \choose n-1} - {n+d-3 \choose n-1}.
$$
We have a natural direct sum decomposition
\begin{equation}\label{L2}
L^2(S^{n-1}) = \bigoplus_{d \ge 0} \mathcal H_{n, d},    
\end{equation}
let $\operatorname{proj}_{n,d}: L^2(S^{n-1}) \rightarrow \mathcal H_{n, d}$ denote the orthogonal projection on the $d$-th component in (\ref{L2}). Given a function $f \in L^2(S^{n-1})$ we write $f^{=d} = \operatorname{proj}_{n,d}(f)$ for brevity.

Fix $t \in [-1,1]$ and let $x, y \in S^{n-1}$ be a pair of points with scalar product $\langle x, y\rangle = t$. We define a bilinear operator $G_t$ on $L^2(S^{n-1})$ as follows. Given functions $f, h \in L^2(S^{n-1})$, let
$$
G_t(f, h) = \int_{SO(n)} f(g x) h(g y) dg,
$$
where $dg$ is the Haar probability measure on the (compact) Lie group $SO(n)$. It is easy to see that this definition does not depend on the choice of points $x$ and $y$.  Informally speaking, $G_t(f, h)$ is the average value of the product $f(x) h(y)$ over all pairs $(x, y) \in S^{n-1} \times S^{n-1}$ with fixed scalar product $t$. In particular, note that $G_1(f, h) = \langle f, h \rangle$. 
For arbitrary $t$, the operator $G_t$ can be diagonalized in the harmonic basis of $L^2(S^{n-1})$ and we have the following, well-known, expansion:
\begin{equation}\label{gegen}
    G_t(f, h) = \sum_{d \ge 0} P_{n, d}(t) \langle f^{=d}, h^{=d}\rangle.
\end{equation}

Here $P_{n, d}(t)$ is the family of {\em Gegenbauer} or {\em ultraspherical polynomials}. They can be uniquely identified by the following two properties:
\begin{itemize}
    \item For any $d\ge 0$, $P_{n,d}$ is a degree $d$ polynomial and we have $P_{n,d}(1) =1$,
    \item For $d\neq d'$ the polynomials $P_{n,d}$ and $P_{n,d'}$ are orthogonal with respect to the measure $(1-t^2)^{\frac{n-3}{2}}dt$ on the interval $[-1,1]$.
\end{itemize}
These polynomials have the following explicit formula (\cite{Abram}, Chapter 22):
\begin{equation}\label{explicit}
    P_{n, d}(t) = \frac{1}{{n+d-3 \choose d}}\sum_{\ell = 0}^{[d/2]} (-1)^{\ell} \frac{\Gamma(d - \ell + \frac{n-2}{2})}{\Gamma(\frac{n-2}{2})\ell! (d-2\ell)!} (2t)^{d-2\ell}.
\end{equation}
For instance, the first three Gegenbauer polynomials are given by
$$
P_{n,0}(t) = 1, ~ P_{n,1}(t) = t, ~ P_{n, 2}(t) = \frac{n}{n-1}t^2 - \frac{1}{n-1}.
$$
We will use some simple bounds on the values of Gegenbauer polynomials at $t=0$. Note that $P_{n,d}(0) = 0$ for odd $d$, since $P_{n,d}(t)$ is an odd function in this case. 
For even $d$ we have 
\begin{equation*}\label{eq:Pexplicit}
P_{n,d}(0) = (-1)^{d/2} \frac{{(n-4)/2 + d/2 \choose d/2}}{{n-3+d \choose d }}.    
\end{equation*}
Using the formula ${a+b \choose b} = \frac{a+b}{b} {a+b-1 \choose b-1}$ a couple of times we get that for $d\ge 2$:
\begin{equation}\label{eq:Precursive}
P_{n,d}(0) = - \frac{d-1}{n-3+d} P_{n,d-2}(0),    
\end{equation}
so for $d\le 6$ we have,
$$
P_{n, 2}(0) = -\frac{1}{n-1}, ~P_{n, 4}(0) = \frac{3}{n^2-1},~P_{n, 6}(0) = -\frac{15}{(n^2-1)(n+3)}.
$$
and for any $d \ge 4$ and $n \ge 2$ we have
\begin{equation}\label{smallP}
    |P_{n, d}(0)| \le |P_{n,4}(0)| \le \frac{3}{n^2-1}.
\end{equation}

For every even $d \ge 2$ and $n\ge 2$ we have:
\begin{equation}\label{eq:Pdecay}
    |P_{n,d}(0)| \le \left(\frac{d}{n}\right)^{d/2}.
\end{equation}
This follows from (\ref{eq:Precursive}) by induction starting from $d=2$. 

Using (\ref{gegen}), we can now give a second proof of Proposition \ref{prop:easy}:
\begin{proof}[Second proof of Proposition \ref{prop:easy}]
    Let $A \subset S^{n-1}$ be a set of density $\alpha$ and let $f = 1_A$. Let $x_1, x_2$ be a uniformly random pair of orthogonal vectors on $S^{n-1}$. Then by definition we have
    \[
    G_0(f, f) = \E f(x_1) f(x_2) = \int f(x_1)  \left( \int f(x_2) d\mu_{S^{n-2}}(x_2)\right) d\mu_{S^{n-1}}(x_1) = \alpha \E_{x \in A} \mu_{S^{n-2}}(A\cap x^\perp)
    \]
    where the second integral is taken over $S^{n-1} \cap x_1^\perp$. On the other hand, by (\ref{gegen}), we have
    \[
    G_0(f, f) = \sum_{d \ge 0} P_{n, d}(0) \|f^{=d}\|_2^2 = \|f^{=0}\|_2^2 - \frac{1}{n-1} \|f^{=2}\|_2^2 + \frac{3}{n^2-1}\|f^{=4}\|_2^2 - \frac{15}{(n^2-1)(n+3)}\|f^{=6}\|_2^2 +\ldots
    \]
    Note that $f^{=0}= \alpha$ and so $\|f^{=0}\|_2^2 = \alpha^2$. On the other hand, by orthogonality we have 
    \[
    \alpha=\|f\|_2^2 = \sum_{d\ge 0} \|f^{=d}\|_2^2 = \alpha^2 + \sum_{d\ge 1}\|f^{=d}\|_2^2.
    \]
    So using $|P_{n,d}(0)| \le \frac{1}{n-1}$ for all $d \ge 2$, we obtain
    \[
    G_0(f, f) \ge \alpha^2 - \frac{1}{n-1} \sum_{d\ge 2}\|f^{=d}\|_2^2 \ge \alpha^2 - \frac{\alpha-\alpha^2}{n-1}.
    \]
    By rearranging, we get the desired lower bound on the expected measure $\mu_{S^{n-2}}(A\cap x^\perp)$.
\end{proof}

\subsection{Hypercontractivity} 
Our argument relies on the classical hypercontractive inequality in Gaussian space proved by Bonami \cite{Bonami}, Gross \cite{Gross} and Beckner \cite{Beckner}. Let $\gamma_n$ denote the standard Gaussian measure on the space $\R^n$. For $\rho \in [0,1]$ define the {\em noise operator} on the space of functions $L^p(\R^n, \gamma_n)$ by
\begin{equation}\label{eq:noise_op}
T_\rho f(x) = \E_{y \sim \gamma_n} \left[ f(\rho x + \sqrt{1-\rho^2} y) \right].    
\end{equation}

\begin{theorem}\label{hyper}
    Let $q \ge p \ge 1$ and let $\rho \le \sqrt{\frac{p-1}{q-1}}$. For any $f \in L^p(\R^n, \gamma_n)$, we have $\|T_\rho f\|_q \le \|f\|_p$.
\end{theorem}

In fact, we only need a simple corollary of Theorem \ref{hyper}: given a degree $d$ harmonic function $f \in \mathcal H_{n,d}$ and $q \ge 2$, we have
\begin{equation}\label{moment}
    \|f\|_{L^q(\R^n, \gamma_n)} \le (q-1)^{d/2} \|f\|_{L^2(\R^n, \gamma_n)}.
\end{equation}
Indeed, this follows from Theorem \ref{hyper} applied to $f$ with $p = 2$ and $\rho = \frac{1}{\sqrt{q-1}}$ and the fact that $T_{\rho} f = \rho^d f$ for any degree $d$ harmonic function $f \in \mathcal H_{n, d}$. For the latter, note that the definition of the noise operator (\ref{eq:noise_op}) can be written as a double integral 
\[
T_\rho f(x) = \int h(r)\int_{S_r(x)} f(y) d\mu_{S_r(x)}(y)  dr
\]
where $S_r(x)$ is the sphere of radius $r$ around $x$ and $h(r)$ is some weight coefficient whose precise form is not relevant to us. Since harmonic functions $f$ satisfy $\int_{S_r(x)} f(y) d\mu_{S_r(x)}(y) = f(x)$, we obtain $T_\rho f(x)  =f(\rho x)$ for any harmonic $f$.
So if $f$ is a harmonic homogeneous degree $d$ polynomial, then $T_\rho f(x) = f(\rho x) = \rho^d f(x)$, as claimed.

For a homogeneous degree $d$ function $f$ the $L^q$-norms in the Gaussian space and on the sphere are related as follows:
\[
\|f\|_{L^q(S^{n-1}, \mu)} =  \left( \frac{\Gamma(\frac{n}{2})}{2^{\frac{dq}{2}} \Gamma(\frac{dq+n}{2})} \right)^{1/q} \|f\|_{L^q(\R^n, \gamma_n)}.
\]
So, converting (\ref{moment}) to the $L_q$-norms on the sphere implies that for any degree $d$ harmonic polynomial $f\in \mathcal H_{n,d}$ and $q\ge 2$:
\begin{equation}\label{eq2}
\|f\|_{L^q(S^{n-1}, \mu)} \le (q-1)^{d/2} e^{\frac{d^2 q}{n}} \|f\|_{L^2(S^{n-1}, \mu)}.    
\end{equation}
An important corollary of this inequality is the so-called {\em Level-$d$ inequality} (see e.g. \cite{ODonnell} for more details in the context of analysis of boolean functions on discrete product spaces):

\begin{prop}\label{leveld}
    Let $n\ge 2$ and $f \in L^2(S^{n-1})$ be a 0-1 valued measurable function. Let $\alpha = \E f$ and suppose that $1/2 \ge \alpha \ge 1/2^n$. Then for any $0 \le d \le \log(1/\alpha)$ and some constant $C_{\ref{leveld}}$ we have
    \begin{equation}\label{ineqd}
    \|f^{=d}\|_2^2 \le \alpha^2 \left( \frac{C_{\ref{leveld}} \log (1/\alpha)}{d} \right)^d.
    \end{equation}
\end{prop}

\begin{proof}
    Let $q \ge 2$ and $q' \in (1,2]$ be dual exponents which we optimize later, then by H{\"o}lder's inequality
    $$
    \|f^{=d}\|_2^2 = \langle f^{=d}, f\rangle \le \|f^{=d}\|_q \|f\|_{q'}.
    $$
    The first term is at most $(q-1)^{d/2} e^{\frac{d^2 q}{n}}\|f^{=d}\|_2$ by (\ref{moment}) and the second term equals $\alpha^{1/q'} = \alpha^{1-1/q}$ since $f$ is a 0-1 valued function with mean $\alpha$. So for any $q \ge 2$ we get
    \[
    \|f^{=d}\|_2^2 \le (q-1)^{d} e^{\frac{2d^2 q}{n}} \alpha^{2-2/q}.
    \]
    Taking $q = \frac{2\log(1/\alpha)}{d}$ we get
    \[
    \|f^{=d}\|_2^2 \le \left(\frac{2\log(1/\alpha)}{d}\right)^{d} e^{\frac{4d \log(1/\alpha)}{n}} \alpha^{2} 2^d \le \left(\frac{4 \log(1/\alpha)}{d}\right)^{d} (1/\alpha)^{\frac{4d}{n}} \alpha^2 \le \left(\frac{64 \log(1/\alpha)}{d}\right)^{d} \alpha^2 
    \]
    where we used the assumption $\alpha \ge 2^{-n}$. This concludes the proof with $C_{\ref{leveld}} = 64$.
\end{proof}

%\subsection{Distribution of quadratic functions} 
We also use (\ref{moment}) to deduce that mean zero quadratic functions are non-positive on a constant fraction of the sphere:

\begin{prop}\label{cor1}
    There is an absolute constant $c_{\ref{cor1}} > 0$ such that for any $n \ge 2$ and a degree 2 harmonic function $f \in \mathcal H_{n, 2}$ we have
    $$
    |\{ x\in S^{n-1}:~ f(x) \le 0 \}| \ge c_{\ref{cor1}}.
    $$
\end{prop}

\begin{proof}
Let $X = (X_1, \ldots, X_n)$ be a sequence of $n$ independent Gaussians $X_i \sim \mathcal{N}(0,1)$ and let $Y = f(X)$. Then the measure of the set of $x\in S^{n-1}$ with $f(x)\le 0$ equals the probability that $Y\le 0$.

By the assumption that $f \in \mathcal H_{n,2}$, we have $\E Y = 0$. A special case of Lemma 5.9 in \cite{Kwan} (see also \cite{Alon}) implies that if $\E Y^4 \le B (\E Y^2)^2$ holds for some constant $B$, then $\Pr[Y \le 0] \ge \frac{1}{5B}$. By (\ref{moment}) applied with $q=4$, we have
$$
\E Y^4 = \|f\|_{L^4(\R^n, \gamma_n)}^4 \le 3^{4} \|f\|_{L^2(\R^n, \gamma_n)}^4 \le 3^4 (\E Y^2)^2, 
$$
and so the result follows with $c_{\ref{cor1}} = \frac{1}{3^4  5} \approx 0.002$.
\end{proof}

\section{Proofs of Lemma \ref{lemperp} and \ref{lem:density_lem}}

\begin{proof}[Proof of Lemma \ref{lemperp}]
    
Now we proceed to the proof of Lemma \ref{lemperp}.
We define $\varepsilon_{\ref{lemperp}} = c_{\ref{cor1}}/2$ and $C_{\ref{lemperp}} = 4\sqrt{2 / c_{\ref{cor1}}}$.
Let $A \subset S^{n-1}$ be a measurable centrally symmetric set of density $\alpha = 1- \varepsilon$ for some $\varepsilon \le \varepsilon_{\ref{lemperp}}$. 
Let $f = 1_A$ be the indicator function of $A$ and consider the expansion 
$$
f = \sum_{d \ge 0} f^{= d},
$$
where $f^{=d} \in \mathcal H_{n,d}$ is a harmonic function of degree $d$. Since $f$ is an even function, $f^{=d} = 0$ for odd $d$. 
%By the Level-$d$ inequality (Lemma \ref{leveld}) applied to the 0-1 valued function $1 - f$ and even $d$, we have
%$$
%\|f^{=d}\|_2^2 = \| (1-f)^{=d} \|_2^2 \le \left( \frac{C \log(1/\varepsilon)}{d}\right)^d \varepsilon^2.
%$$

By Proposition \ref{cor1} applied to $f^{=2}$, we have
\[
|\{x \in S^{n-1}: f^{=2}(x) \le 0\}| \ge c_{\ref{cor1}}
\]
and so if we let $B = \{x \in A:~ f^{=2}(x) \le 0\}$ then we have 
\[
|B| \ge c_{\ref{cor1}} - |S^{n-1}\setminus A| \ge c_1 -\varepsilon_{\ref{lemperp}} \ge c_{\ref{cor1}}/2,
\]
if we take $\varepsilon_{\ref{lemperp}} \le c_{\ref{cor1}}/2$.
Let $h = 1_B$ and denote $\beta = |B|$. Observe that then by definition
\begin{equation}\label{deg2}
\langle f^{=2}, h^{=2} \rangle = \langle f^{=2}, h \rangle = \IS f^{=2}(x) 1_B(x) d\mu(x) \le 0,    
\end{equation}
so we can expand
\begin{align*}
G_0(f, h) = \sum_{d\ge 0} P_{n,d}(0) \langle f^{=d}, h^{=d} \rangle = |A| |B| - \frac{1}{n-1} \langle f^{=2}, h^{=2} \rangle + \frac{3}{n^2-1}  \langle f^{=4}, h^{=4} \rangle + \ldots \ge \\
\ge \alpha\beta - \frac{3}{n^2-1} \sum_{d\ge 4, \text{ even}} |\langle f^{=d}, h^{=d} \rangle|,
\end{align*}
where we used (\ref{smallP}) in the end. 
By the Cauchy--Schwarz inequality and the fact that $(1-f)^{=d} = -f^{=d}$ for $d \neq 0$, we have
\[
\sum_{d\ge 4, \text{ even}} |\langle f^{=d}, h^{=d} \rangle| \le \left( \sum_{d\ge 4, \text{ even}} \|(1-f)^{=d}\|_2^2 \right)^{1/2} \left( \sum_{d\ge 4, \text{ even}} \|g^{=d}\|_2^2 \right)^{1/2} \le \|1-f\|_2 \|g\|_2 \le \varepsilon^{1/2} \beta^{1/2}.
\]
and thus,
$$
G_0(f, h) \ge \alpha \beta - \frac{3 \varepsilon^{1/2} \beta^{1/2}}{n^2-1}.
$$

On the other hand, we can write
$$
G_0(f, h) = \IS 1_B(x) \mu_{S^{n-2}}(A \cap x^\perp) d\mu(x) = \beta~ \E_{x\in B}~ \mu_{S^{n-2}}(A \cap x^\perp),
$$
where $x^\perp$ is the hyperplane in $\R^n$ orthogonal to $x$ (note that the measure of the intersection $A \cap x^\perp$ is well-defined for almost every $x \in A$). Thus, we conclude that there exists a point $x \in B\subset A$ such that 
$$
\mu_{S^{n-2}}(A \cap x^\perp) \ge \alpha -\frac{3 \varepsilon^{1/2} \beta^{-1/2}}{n^2-1} \ge \alpha- \frac{4 \varepsilon^{1/2} \beta^{-1/2}}{n^2}
$$
using $n\ge 2$.
Recall that $\beta = |B| \ge c_{\ref{cor1}}/2$ so $\beta^{-1}$ is bounded by an absolute constant. So we get the lower bound claimed in Lemma \ref{lemperp} with $C_{\ref{lemperp}} = 4\sqrt{2 / c_{\ref{cor1}}}$.

\end{proof}

\begin{proof}[Proof of Lemma \ref{lem:density_lem}]
    Let $A \subset S^{n-1}$ be a set of density $\alpha$. 
    We will put $C_{\ref{lem:density_lem}} = 2C_{\ref{leveld}}^2+1$. 
    
    Let $f = 1_A$ then we have 
    \begin{equation}\label{eq:expandG}
    \E_{x \in A} \mu_{S^{n-2}}(A \cap x^\perp) = \alpha^{-1} G_0(f, f) = \alpha^{-1} \sum_{d = 0}^{\infty} P_{n,d}(0) \|f^{=d}\|_2^2,    
    \end{equation}
    so it is sufficient to lower bound the right hand side. 

    %We split the proof in two cases depending on how large is the density $\alpha$. 
    %First, consider the case when $\alpha \in (0, e^{-2})$.
    Let $d_0 := 2[\log(1/\alpha)/2]$, note that $d_0\ge 2$ by the assumption that $\alpha \le e^{-2}$.
    By the Level $d$ inequality, we have for all $d \le d_0$:
    \[
    \|f^{=d}\|_2^2 \le \alpha^2 \left ( \frac{C_{\ref{leveld}}\log (1/\alpha)}{d} \right)^d.
    \]
    So by (\ref{eq:Pdecay}) we have
    \[
    \sum_{d=2}^{d_0} |P_{n,d}(0)| \|f^{=d}\|_2^2 \le \alpha^2 \sum_{d=2}^{d_0} (d / n)^{d/2} \left ( \frac{C_{\ref{leveld}}\log (1/\alpha)}{d} \right)^d \le \alpha^2 \sum_{d=2}^{d_0} \left ( \frac{C_{\ref{leveld}}\log (1/\alpha)}{\sqrt{n}} \right)^d 
    \]
    So provided that $\log (1/\alpha) \le \frac{\sqrt{n}}{2 C_{\ref{leveld}}}$, summing the geometric series gives
    \[
    \sum_{d=2}^{d_0} |P_{n,d}(0)| \|f^{=d}\|_2^2 \le \alpha^2 \frac{2 C^2_{\ref{leveld}}\log^2 (1/\alpha)}{n}.
    \]
    Note that if $\log (1/\alpha) > \frac{\sqrt{n}}{2 C_{\ref{leveld}}}$ then the bound in Lemma \ref{lem:density_lem} is trivially true.
    
    For $d> d_0$ we can use (\ref{eq:Pdecay}) to estimate
    \[
    \sum_{d=d_0+2}^\infty |P_{n,d}(0)| \|f^{=d}\|_2^2 \le (d_0 / n)^{d_0/2}  \sum_{d=d_0+2}^\infty \|f^{=d}\|_2^2 \le (d_0 / n)^{d_0/2} \|f\|_2^2 = (d_0 / n)^{d_0/2}\alpha
    \]

    So plugging these estimates in (\ref{eq:expandG}) we get
    \[
    \E_{x \in A} \mu_{S^{n-2}}(A \cap x^\perp) \ge \alpha- \alpha \frac{2 C^2_{\ref{leveld}}\log^2 (1/\alpha)}{n} - \alpha(d_0/n)^{d_0/2} \ge \alpha\left(1- \frac{(2C_{\ref{leveld}}^2+1) \log^2(1/\alpha)}{n} \right)
    \]
    as desired.
\end{proof}

\section{Proofs of Theorem \ref{main} and Theorem \ref{prop:other_simplices}}

\begin{proof}[Proof of Theorem \ref{main}]

We prove the result using induction on $n$.
Let $\varepsilon = \min\{\frac{\varepsilon_{{\ref{lemperp}}}}{2}, \frac{1}{32C^2_{\ref{lemperp}}}\}$, we are going to prove the result for $c_0 = 1-\varepsilon$.
Denote $n_0 = [ \frac{1}{2\varepsilon} ]$. For $n \le n_0$ we have the easy bound (see Introduction for the proof)
\begin{equation}\label{eq:easy-bound}
m_{S^{n-1}}(\Delta_{n,0}) \le 1 - 1/n \le 1-2\varepsilon 
\end{equation}
so we may assume that $n>n_0$.

Let $A \subset S^n$ be a set of density at least $1-\varepsilon$. Note that if $A\cup (-A)$ contains $n$ pairwise orthogonal vectors then so does the set $A$. So without loss of generality, we may assume that $A$ is centrally symmetric.
Let $m = n-n_0$ and use Lemma \ref{lemperp} to construct a sequence of pairwise orthogonal points $x_1, \ldots, x_m \in A$ so that if we denote 
\[
\mu_{S^{n-1-j}}(A \cap x_1^\perp \cap \ldots \cap x_j^\perp) = 1-\varepsilon_{j} 
\]
then we have for every $j = 0, \ldots, m-1$:
\[
1-\varepsilon_{j+1} \ge 1 - \varepsilon_{j} - \frac{C_{\ref{lemperp}} \varepsilon_{j}^{1/2}}{(n-j)^2},
\]
where we also denote $\varepsilon_0 = \varepsilon$. 
First, we claim that $\varepsilon_j \le 2\varepsilon$ for every $j \le m$ (note that this would imply $\varepsilon_j \le \varepsilon_{\ref{lemperp}}$ and so the application of Lemma \ref{lemperp} at step $j$ is justified). Indeed, if we already know that $\varepsilon_i \le 2\varepsilon$ for $i = 0, \ldots, j-1$ then we get
\[
\varepsilon_{j} \le \varepsilon + C_{\ref{lemperp}} \sum_{i=0}^{j-1} \frac{\varepsilon_i^{1/2}}{(n-i)^2} \le \varepsilon + C_{\ref{lemperp}} (2\varepsilon)^{1/2} \sum_{i=0}^{j-1} \frac{1}{(n-i)^2} \le \varepsilon + \frac{C_{\ref{lemperp}} (2\varepsilon)^{1/2}}{n-j} \le \varepsilon + \frac{C_{\ref{lemperp}} (2\varepsilon)^{1/2}}{n_0}
\]
So we get $\varepsilon_j \le 2\varepsilon$ provided that $n_0 \ge C_{\ref{lemperp}} \sqrt{2/\varepsilon}$ holds. This condition is satisfied by our choice of $\varepsilon$ and $n_0$. 

Taking $j = n-n_0$, we conclude that the intersection $B = A\cap x_1^\perp \cap \ldots \cap x_{n-n_0}^{\perp}$ satisfies
\[
\mu(B) \ge 1-2\varepsilon > 1- \frac{1}{n_0}.
\]
By (\ref{eq:easy-bound}), the set $B$ contains $n_0$ pairwise orthogonal vectors $y_1,\ldots, y_{n_0}$. Combining them with the earlier constructed sequence of points $x_1, \ldots, x_{n-n_0} \in A$ gives us the desired configuration of points. This completes the proof of Theorem \ref{main} using Lemma \ref{lemperp}.

\end{proof}

\begin{proof}[Proof of Theorem \ref{prop:other_simplices}]

First we observe that the range $k \in [c n, n]$ follows directly from Theorem \ref{main}: if $A \subset S^{n-1}$ does not contain $k$ pairwise orthogonal vectors then by Theorem \ref{main}, we have $\mu(A) \le c_0$. So if $k\ge cn$ then we get $\mu(A) \le \exp(-c'_1 n/k)$ with $c'_1 = c(1-c_0)$. So we may assume that $k \le cn$ for any fixed constant $c>0$. We will choose $c = \frac{1}{32 C_{\ref{lem:density_lem}}}$. Note that this in particular implies that $n \ge k/c \ge 64 C_{\ref{lem:density_lem}}$.

Let $A \subset S^{n-1}$ be a set of density $\alpha$ such that 
\begin{equation}\label{eq:condition-alpha}
\alpha \ge \exp\left(- \min\left( (1/16 C_{\ref{lem:density_lem}}) n/k, \sqrt{n / 16 C_{\ref{lem:density_lem}}} \right)\right)    
\end{equation}
By shrinking the set $A$ if necessary we may also assume that $\mu(A) \le e^{-2}$ holds (note that the restrictions on $n$ and $k$ guarantee that this does not conflict with (\ref{eq:condition-alpha})).

Now we use Lemma \ref{lem:density_lem} iteratively to obtain bounds on sets avoiding $\Delta_{k, 0}$. Indeed, let $x_1, \ldots, x_{k-1} \in A$ be a sequence obtained by $k-1$ applications of Lemma \ref{lem:density_lem} and for $j=0, \ldots, k-1$ let 
\[
A_j = A \cap x_1^\perp \cap \ldots \cap x_j^\perp.
\]
Write $\alpha_j = \mu_{S^{n-1-j}}(A_j)$, then Lemma \ref{lem:density_lem} gives for $0\le j \le k-2$:
\[
\alpha_{j+1} \ge \alpha_j \left(1- \frac{C_{\ref{lem:density_lem}} \log^2(1/\alpha_j)}{n-j} \right).
\]
We claim that $\alpha_j \ge \alpha_0^2$ for all $j =1, \ldots, k-1$. Indeed, suppose that $\alpha_i \ge \alpha_0^2$ holds for all $i \le j-1$ then we get
\[
\alpha_{j} \ge \alpha_{j-1} \left(1- \frac{C_{\ref{lem:density_lem}} \log^2(1/\alpha_{j-1})}{n-j+1} \right) \ge \alpha_{j-1} \left ( 1- \frac{8 C_{\ref{lem:density_lem}}\log^2 (1/\alpha_0)}{n} \right) \ge \ldots \ge \alpha_0 \left(1- \frac{8 C_{\ref{lem:density_lem}} \log^2(1/\alpha_0)}{n} \right)^{j}
\]
where we used $n-j+1 \ge n/2$. So provided that $\log (1/\alpha_0) \le \sqrt{n/16 C_{\ref{lem:density_lem}}}$, we get
\[
\alpha_{j} \ge \alpha_0 \exp(- 16 C_{\ref{lem:density_lem}} \log^2(1/\alpha_0) j/ n)
\]

So if $\alpha_0 \ge \exp(- \frac{n}{ 16 C_{\ref{lem:density_lem}} k})$ then this implies that $\alpha_{j+1} \ge \alpha_0^2$ as desired. Both conditions on $\alpha_0$ are indeed satisfied by (\ref{eq:condition-alpha}). In particular we get $\alpha_{k-1} \ge \alpha_0^2 > 0$ and we obtain that $A$ contains $k$ pairwise orthogonal vectors.

So for $k \le c n$, we conclude that
\[
m_{S^{n-1}}(\Delta_{k,0}) \le \exp\left(- c_1'' \min\left( n/k, \sqrt{n} \right)\right)   
\]
with $c_1'' = 1/16C_{\ref{lem:density_lem}}$. So combining with the range $k \in [cn, n]$ we obtain Theorem \ref{prop:other_simplices} with $c_1 = \min(c_1', c_1'') = \frac{1-c_0}{32 C_{\ref{lem:density_lem}}}$, completing the proof.

\end{proof}

\bibliographystyle{amsplain0.bst}
\bibliography{main}

\end{document}